\documentclass{amsart}[12pt]
\usepackage{amsmath, amstext, amssymb, amsopn, amsthm, amscd, amsfonts, 
mathptmx, epsfig, times}

\newcommand{\R}{\mathbb R}

\newcommand{\HP}{\mathbb H}

\newcommand{\N}{\mathbb N}

\newtheorem{theo}{Theorem}
\newtheorem{lemm}{Lemma}
\newtheorem{prop}{Proposition}

\newtheorem*{echyp}{Ehrenpreis Conjecture (Hyperbolic)}
\newtheorem*{ecconf}{Ehrenpreis Conjecture (Conformal)}
\newtheorem*{ecsol}{Ehrenpreis Conjecture (Solenoidal)}
\newtheorem*{l1ec}{$L^{1}$ Ehrenpreis Conjecture}

\theoremstyle{definition}

\theoremstyle{remark}

\newtheorem{note}{Note}

\newtheorem{clai}{Claim}

\title{The $L^{1}$ Ehrenpreis Conjecture}
\author{T.M. Gendron}
\address{Instituto de Matem\'{a}ticas -- Unidad Cuernavaca, Universidad
Nacional Autonoma de M\'{e}xico, Av. Universidad S/N, C.P. 62210
Cuernavaca, Morelos, M\'{E}XICO}
\email{tim@matcuer.unam.mx}
\date{12 June 2005}
\begin{document}
\vspace{2cm} \maketitle
\begin{abstract} Let $\widehat{\Sigma}$ be the algebraic universal
cover of a closed surface of genus $g>1$, 
$\mathcal{T}(\widehat{\Sigma})$ its Teichm\"{u}ller space,
${\sf M}(\Sigma ,\ell )$
the group of mapping classes stabilizing a fixed leaf $\ell$.  The
$L^{1}$ Ehrenpreis conjecture asserts that ${\sf M}(\Sigma ,\ell )$
acts on $\mathcal{T}(\widehat{\Sigma})$ with dense orbits in the 
$L^{1}$ topology (the topology coming from the $L^{1}$ norm on quadratic differentials).
We give a proof of this weaker version of the Ehrenpreis conjecture,
announced first in \cite{Ge1}.  
\end{abstract}

\section{Introduction}

Let $\Sigma, \Sigma'$ be closed Riemann surfaces of genus greater than 1.  
The most succinct formulation of the Ehrenpreis conjecture ({\sf EC}) 
uses the fact that $\Sigma, \Sigma'$ may be regarded as riemannian manifolds
with metrics of curvature $-1$.
While it is an elementary fact that the riemannian universal covers 
$\widetilde{\Sigma}, \widetilde{\Sigma}'$ are isometric to $\HP^{2}$, the 
{\sf EC}
asserts a similar, asymptotic phenomenon
for the family of finite riemannian covers: 
\begin{echyp}  For any $\epsilon >0$, there are finite degree 
isometric covers $Z\rightarrow \Sigma$, $Z'\rightarrow \Sigma'$ whose
total spaces are $(1+\epsilon )$-quasiisometric.
\end{echyp}
In \cite{Ge1}, we announced the solution of an
$L^{1}$-version of this conjecture.  In this paper, we provide the proof.

The traditional or conformal version of the {\sf EC} \cite{Eh} can
be described in terms of
Teichm\"{u}ller theory.  
Let $\Sigma$ be a fixed compact surface of genus greater than 1, 
$\mathcal{T}(\Sigma )$ its Teichm\"{u}ller space, 
$d_{\mathcal{T}(\Sigma )}$ the Teichm\"{u}ller metric.  Then given
$\mu ,\nu\in \mathcal{T}(\Sigma )$, the conformal
version of the {\sf EC} states 

\begin{ecconf} For any $\epsilon >0$, there exists a 
surface $Z$ and finite covers $\rho ,\sigma :Z\rightarrow
\Sigma $ such that
\[  d_{\mathcal{T}(Z )}(\rho^{\ast}\mu , \sigma^{\ast}\nu )\; <\; \epsilon .\]
\end{ecconf}

In the above, $\rho^{\ast}, \sigma^{\ast} :\mathcal{T}(\Sigma )\hookrightarrow\mathcal{T}(Z )$
are the isometric inclusions induced by $\rho ,\sigma$.  We remark
that this version of the {\sf EC} makes sense in genus 1, where it is
not difficult to verify \cite{Ge1}.

However, it is the genus independent or solenoidal version of the {\sf EC} that will be most
important for us.  
Let  $\widehat{\Sigma}$
be the algebraic universal cover of the closed surface $\Sigma$, by definition
the inverse limit of the total spaces of finite covers $\rho :Z\rightarrow \Sigma$
(one cover for each homotopy class of cover).  The algebraic universal cover
is a surface solenoid (a surface lamination with Cantor transversals), and as
such has a Teichm\"{u}ller space $\mathcal{T}(\widehat{\Sigma})$
of marked conformal structures \cite{Su}.  The mapping class group 
${\sf M}(\widehat{\Sigma},\ell )$
of homotopy classes of homeomorphisms of $\widehat{\Sigma}$ fixing a base leaf $\ell$
may be identified with the group of homotopy classes of lifts of correspondences
$\sigma\circ\rho^{-1}:\Sigma\rightarrow\Sigma$.  

\begin{ecsol}
${\sf M}(\widehat{\Sigma}, \ell )$ acts on $\mathcal{T}(\widehat{\Sigma})$
with dense orbits.
\end{ecsol}
The proof that all of these versions of the {\sf EC} are equivalent can be found
in \cite{Ge1}.

The genus independent version of {\sf EC} says
that the "universal moduli space"
\[ \mathcal{M}(\widehat{\Sigma}) \;\; =\;\;
\mathcal{T}(\widehat{\Sigma})/{\sf M}(\widehat{\Sigma}, \ell ),\] although
an uncountable set, has the topology of a point (the discrete topology).  
It has the virtue of
giving a certain explanation of the ``moduli-rigidity gap'' that separates the
theory of compact hyperbolic surfaces from that of compact hyperbolic manifolds in
dimension three or greater.  
From a practical point of view, working with $\widehat{\Sigma}$
allows us to regard Teichm\"{u}ller theory of closed hyperbolic surfaces as
concerning complex structures on a single topological type (as in the
genus 1 case).  In this way, we may isolate geometric properties of 
closed Riemann surfaces of hyperbolic type
that do not depend on genera. 

In this paper, we shall formulate and prove an $L^{1}$ version of the {\sf EC}.
On $\mathcal{T}(\widehat{\Sigma})$ there are -- in addition to the Teichm\"{u}ller
metric -- three other metrics coming from the $L^{1}$, the $L^{\infty}$ and
the $L^{2}$ structures on the cotangent bundle of $\mathcal{T}(\widehat{\Sigma})$.
The $L^{1}$ version of the {\sf EC} is obtained by asking that 
${\sf M}(\widehat{\Sigma},\ell )$
act densely on $\mathcal{T}(\widehat{\Sigma})$ with regard to the $L^{1}$ geometry.

The proof of the $L^{1}$ {\sf EC} is as follows.  A dense set of pairs
$\hat{\mu} ,\hat{\nu} \in \mathcal{T}(\widehat{\Sigma})$ (dense in the Teichm\"{u}ller
geometry) lie along
the axis $A$ of a pseudo-Anosov homeomorphism $\hat{\Phi} :\widehat{\Sigma}\rightarrow
\widehat{\Sigma}$.  $A$ is a Teichm\"{u}ller geodesic and the
action of $\hat{\Phi}$ on $\mathcal{T}(\widehat{\Sigma})$ stabilizes
$A$, translating points along $A$ a distance of $\frac{1}{2}\log\lambda$,
where $\lambda$ is the entropy of $\hat{\Phi}$.
Given $n\in\N$, by an $L^{1}$ $n$th root of $\Phi$ we mean a sequence
of pseudo Anosov homeomorphisms $\{ \hat{\Psi}_{m}\}$ in which 
\begin{itemize}
\item The entropies $\lambda_{m}$ of $\hat{\Psi}_{m}$
converge to $\lambda^{1/n}$.
\item The axes $A_{m}$ of $\hat{\Psi}_{m}$
converge to $A$ in the $L^{1}$ Hausdorff topology.
\end{itemize}
We shall show that $L^{1}$ $n$th roots exist for every (lifted) pseudo Anosov
homeomorphism of $\widehat{\Sigma}$.  This will then imply
the $L^{1}$ version of {\sf EC}.

\vspace{3mm}

\noindent {\bf Acknowledgements:}  This work benefited greatly from
conversations with Dennis Sullivan and Yair Minksy.

\section{Topology of the Algebraic Universal Cover}

Let $\Sigma$ be a fixed compact surface of genus at least two.  
We describe in this section the topology
of the algebraic universal cover $\widehat{\Sigma}$.  Unless otherwise
noted, all proofs of statements in this section can be found in \cite{Ge1}.

Let $\pi =\pi_{1}\Sigma$.  For every finite index normal subgroup $H<\pi$, 
choose a pointed cover $\rho :(Z,x_{Z})\rightarrow(\Sigma ,x)$
for which $\rho_{\ast}\pi_{1}Z=H$.  By adding to this collection of covers all covers
$\tau:Z\rightarrow Z'$ between total spaces for which $\rho '\circ\tau =\rho$,
we obtain an inverse system of surfaces.  The limit of this system $\widehat{\Sigma}$
is called the {\em algebraic universal cover} of $\Sigma$, a compact topological
space. Its topological
type is independent of the choice of covers.

If we denote by $\hat{\pi}$ the profinite completion of $\pi$, then
$\widehat{\Sigma}$ is homeomorphic to the quotient
\begin{equation}\label{suspension} 
\left( \widetilde{\Sigma}\times\hat{\pi} \right)\Big/ \pi , 
\end{equation}
where $\pi$ acts diagonally, and so has the structure of a surface solenoid:
a surface lamination whose model transversals are Cantor sets.  It
also follows from (\ref{suspension}) that $\widehat{\Sigma}$ is connected, 
its path components are its leaves, and each leaf is
homeomorphic to $\R^{2}$ and dense in $\widehat{\Sigma}$.  The point
$\hat{x}=(x_{Z})\in\widehat{\Sigma}$ -- defined by the string of basepoints of the surfaces in 
the defining system -- is contained in a leaf $\ell$ which
we call the base leaf.  The Haar measure on  $\hat{\pi}$ induces
a transverse invariant measure $\eta$ on $\widehat{\Sigma}$
which gives measure $1/\deg Z$ to the fibers of the natural projection
$\widehat{\Sigma}\rightarrow Z$.

Any pointed finite cover $\sigma :(Y,y)\rightarrow (\Sigma ,x)$ 
lifts to a base leaf preserving
homeomorphism $\hat{\sigma}:\widehat{Y}\rightarrow\widehat{\Sigma}$,
where $\widehat{Y}$ is the algebraic universal cover of $Y$. 
If $\rho :(Y,y)\rightarrow (\Sigma ,x)$  is another such cover, the correspondence
$\sigma\circ\rho^{-1}$ lifts to the homeomorphism $\hat{\sigma}\circ\hat{\rho}^{-1}$
of $\widehat{\Sigma}$ preserving $\ell$.  Let ${\sf M}(\widehat{\Sigma},\ell )$ denote the
the group of homotopy classes of orientation preserving homeomorphisms
of $\widehat{\Sigma}$ preserving $\ell$.  See \cite{Ge2}, \cite{Od}
for a proof of the following

\begin{theo}  Every class $[h]\in {\sf M}(\widehat{\Sigma},\ell )$ contains an element of
the form $\hat{\sigma}\circ\hat{\rho}^{-1}$.
\end{theo}

\section{Measured Laminations on $\widehat{\Sigma}$}

We begin by recalling a few facts about measured laminations on $\Sigma$,
a closed Riemann surface of genus $g>1$.  See \cite{Ca-Bl}, \cite{Fa}, \cite{Ha}, \cite{Th1} for further discussion.
A {\em measured lamination} $\mathfrak{f}$ on $\Sigma$ is a closed 1-dimensional lamination
smoothly embedded in $\Sigma$ and possessing a transverse invariant measure 
$m_{\mathfrak{f}}$. Two measured laminations are {\em equivalent} if they
are isotopic through an isotopy taking one measure to the other. The set of equivalence classes of measured laminations
is denoted $\mathcal{ML}(\Sigma )$.  
Let $\mathcal{C}(\Sigma )$ denote the set of isotopy classes
of simple closed curves in $\Sigma$.  Given $\mathfrak{f}\in
\mathcal{ML}(\Sigma )$ and $c\in \mathcal{C}(\Sigma )$, the 
{\em intersection pairing}
is defined
\[  \text{\bf I}( \mathfrak{f},c)\;\; =\;\; 
\inf\int_{\mathfrak{f}\cap c}dm_{\mathfrak{f}} ,  \]
where the infimum is taken over representatives of 
the classes of $\mathfrak{f}$ and $c$.
The intersection topology on $\mathcal{ML}(\Sigma )$ is the weak topology
with respect to the intersection pairing.  
The space of projective classes of measured laminations is denoted
$\mathcal{PL}(\Sigma )$ and is homeomorphic to a sphere of dimension
$6g-7$.  We have 
$\mathcal{C}(\Sigma )\subset \mathcal{PL}(\Sigma )$
with dense image.
The intersection pairing extends to a map
$\mathcal{ML}(\Sigma )\times \mathcal{ML}(\Sigma )\rightarrow\R$ via
the formula \cite{Ke}
\[ \text{\bf I}( \mathfrak{f},\mathfrak{g})\;\; =\;\;
 \inf\int_{\mathfrak{f}\cap \mathfrak{g}}dm_{\mathfrak{f}}\otimes dm_{\mathfrak{g}} .\]

A word is in order here regarding the allied concept of a {\em measured foliation},
a singular foliation $\mathcal{F}$ 
of $\Sigma$ equipped with a transverse invariant measure: these
typically arise as trajectories of holomorphic quadratic differentials on
$\Sigma$ \cite{St}.  Two measured foliations
are equivalent if after a finite number of Whitehead moves are applied to their
singular leaves, they are isotopic through an isotopy taking
one measure to the other.  There is a 
bijective correpondence between classes of measured foliations and
classes of measured laminations \cite{Ha}. For example, to obtain
a measured lamination starting
with a measured foliation $\mathcal{F}$, one chooses a nonsingular leaf
from each minimal component of $\mathcal{F}$, 
pulls each such leaf geodesic (with respect to the hyperbolic
metric of $\Sigma$) and completes
the resulting space.  Our default will be to work with measured laminations,
and -- with the exception of the proof of Theorem~\ref{intconvimpl1}, where we
revert back to measured foliations -- whenever a measured foliation happens to arise, we will assume
it has been converted into its associated measured lamination.

A (homotopy class of) homeomorphism  $\Phi:\Sigma\rightarrow\Sigma$ 
induces a homeomorphism
of $\mathcal{ML}(\Sigma )$ via pullback of measures, in particular
inducing a homeomorphism of $\mathcal{PL}(\Sigma )$.
According to the classification of surface diffeomorphisms, 
\cite{Ca-Bl}, \cite{Fa}, \cite{Th},
$\Phi$ is called {\em pseudo Anosov} if
its induced action on $\mathcal{PL}(\Sigma )$
fixes precisely two classes $[\mathfrak{f}^{u}]$ and $[\mathfrak{f}^{s}]$.
If $\lambda$ is the entropy of $\Phi$, then $\lambda >1$; and if
$\mathfrak{f}^{u}\in [\mathfrak{f}^{u}]$ ($\mathfrak{f}^{s}\in [\mathfrak{f}^{s}]$) is a representative in
$\mathcal{ML}(\Sigma )$, then there is a representative diffeomorphism
in the class of $\Phi$ (also denoted $\Phi$) such that 
$\Phi (\mathfrak{f}^{u} )=\lambda \mathfrak{f}^{u}$ 
($\Phi (\mathfrak{f}^{s})=\lambda^{-1} \mathfrak{f}^{s}$).

In this paper, we will be interested in the following class of
pseudo Anosov homeomorphisms.  Let $\mathcal{C}$, $\mathcal{D}$ be 
families of pairwise nonisotopic simple closed curves, for which
elements of $\mathcal{C}$ intersect minimally in their isotopy classes
with elements of $\mathcal{D}$, and for 
which $\Sigma\setminus (\mathcal{C}\cup\mathcal{D})$ consists of
a union of disks (the families are then said to be {\em filling}).  
For $c\in\mathcal{C}$,
$d\in\mathcal{D}$, let $F_{c}$ resp.\ $G_{d}$ denote the right Dehn twist
about $c$ resp.\ $d$.  Then a homeomorphism of the form
\[ \Phi\;\; =\;\; G^{-N_{k}}_{d_{k}}\circ\dots\circ G^{-N_{1}}_{d_{1}}\circ
F^{M_{j}}_{c_{j}}\circ\dots\circ F^{M_{1}}_{c_{1}} ,\]
where the exponents $M_{1},\dots , M_{j},N_{1},\dots , N_{k}$ are positive
and
where all curves in $\mathcal{C}$, $\mathcal{D}$ occur,
is pseudo Anosov, \cite{Pe2}, \cite{Th1}.  We call these pseudo Anosovs of
{\em Thurston-Penner type}.

We now extend the above considerations to $\widehat{\Sigma}$.  
A measured lamination $\hat{\mathfrak{f}}$ on $\widehat{\Sigma}$ is a collection of
measured laminations $\{ \mathfrak{f}_{\ell}\}$, one on
each leaf $\ell$ of $\widehat{\Sigma}$, which have the same
transversal model ${\sf T}$ and which
vary in the following way with respect to the transversals of $\widehat{\Sigma}$.
Let $\mathcal{O}\approx D\times \widehat{T}$ be a flowbox for $\widehat{\Sigma}$,
such that $\mathfrak{f}_{t}:=
\hat{\mathfrak{f}}|_{D\times\{ t\}}$ is a flowbox for $ \mathfrak{f}_{\ell}$
if $D\times\{ t\}\subset\ell$.
Then
\begin{enumerate}
\item  The family of flowboxes $\mathfrak{f}_{t}$ varies continuously in $t$. 
(Thus $\hat{\mathfrak{f}}$
gives rise to a smooth 1-dimensional sublamination of $\widehat{\Sigma}$ 
with transversal models $\widehat{T}\times {\sf T}$.)
\item Given any continuous family of test transversals ${\sf T}_{t}\subset
D\times \{ t\}$, the function obtained by pairing with
the measures of $\hat{\mathfrak{f}}$ is continuous in $t$.
\end{enumerate}

The space of equivalence classes of measured laminations is
denoted $\mathcal{ML}(\widehat{\Sigma})$. 
If $c$ is a simple closed curve occurring in some surface $Z$ in the
defining system of $\widehat{\Sigma}$, its preimage $\hat{c}$ in
$\widehat{\Sigma}$ is a $1$-dimensional solenoid which we 
call a {\em simple closed solenoid} (abusively, since such a $\hat{c}$
always has infinitely
many connected components).
The set of such is denoted $\mathcal{C}(\widehat{\Sigma})$.

The intersection pairing $\text{\bf I}(\hat{\mathfrak{f}},\hat{c})$
between a measured lamination and a simple closed solenoid is
defined using the transverse invariant measure $\eta$,
\[  \text{\bf I}(\hat{\mathfrak{f}},\hat{c})\;\; =\;\;
\int_{ \hat{\mathfrak{f}}\cap\hat{c}}dm_{\hat{\mathfrak{f}}}\; d\eta . \]
We equip $\mathcal{ML}(\widehat{\Sigma})$ with
the resulting weak topology.  With its induced
topology, $\mathcal{PL}(\widehat{\Sigma})$ is precompact (being
essentially a space of probability measures), but owing to its infinite-dimensionality,
$\mathcal{PL}(\widehat{\Sigma})$ is 
{\em not} compact.
The simple closed solenoids $\mathcal{C}(\widehat{\Sigma})$
are dense in $\mathcal{PL}(\widehat{\Sigma})$.

If $Z$ is a surface occurring in the defining limit of $\widehat{\Sigma}$,
then we may pullback measured laminations on $Z$ to measured laminations
on $\widehat{\Sigma}$.  The result is a direct system of
inclusions  $\mathcal{ML}(Z)\hookrightarrow 
\mathcal{ML}(\widehat{\Sigma})$ whose limit
$\mathcal{ML}^{\infty}\; :=\;\lim_{\longrightarrow}\mathcal{ML}(Z)$ 
has dense image in $\mathcal{ML}(\widehat{\Sigma})$.  In this paper we will
work exclusively with $\mathcal{ML}^{\infty}$.

A pseudo Anosov diffeomorphism $\Phi :Z\rightarrow Z$ lifts to a diffeomorphism
of $\widehat{\Sigma}$ fixing precisely the lifted projective classes
$[\hat{\mathfrak{f}}^{s}]$ and $[\hat{\mathfrak{f}}^{u}]$.  In what
follows, the terminology 
``pseudo Anosov homeomorphism of $\widehat{\Sigma}$'' will always mean
such a lift.

\section{Train Tracks}

We recall first some facts about train tracks on $\Sigma$: details
may be found in \cite{Pe1}.
Let $\tau\subset\Sigma$ be a 
smooth 1-dimensional branched manifold: thus $\tau$ is a 1-dimensional
CW-complex in which the interiors of edges are smooth curves, and
the field of tangent lines ${\rm T}_{x}\tau$, $x\in\tau\setminus \{\text{vertices}\}$,
extends to a continuous line field on $\tau$.
We say that $\tau$ is a {\em train track} 
if it satisfies the following additional properties:
\begin{enumerate}
\item  The valency of any vertex is at least 3, except for simple closed curve
components, which have a single vertex of valence 2.
\item  If $D(S)$ is the double of a component $S\subset\Sigma\setminus\tau$,
then the Euler
characteristic of $D(S)$ is negative.
\end{enumerate}
We shall follow the custom of referring to the vertices of a train track
as {\em switches}.  
 
 A {\em bigon track} is a smooth 1-dimensional branched manifold $\tau\subset\Sigma$ 
 satisfying
 item (1) and which satisfies (2) after collapsing
 bigon complementary regions to curves.  Bigon tracks arise naturally
 from a pair $\mathcal{C}$, $\mathcal{D}$ of transverse filling curves, 
 by turning each intersection of a $\mathcal{C}$ curve with a $\mathcal{D}$
 curve into a pair of 3-valent vertices as in
 Figure 1.  Since such bigon tracks will be the only ones appearing in this
 article, we will assume from now on that all switches in
 bigon tracks have valency no more than three.
 
 \begin{figure}[htbp]
\centering \epsfig{file=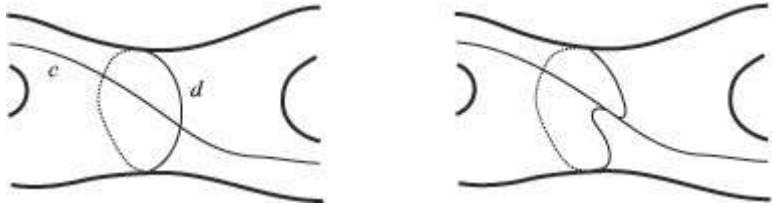, width=4in} \caption{Creating a bigon track
from a filling pair of curves.}
\end{figure}
 
 Denote by $E$ the set of edges of the bigon track $\tau$.  
 In a small disk neighborhood of a switch
 $v$, the ends of edges incident to $v$ may be divided into two classes, which
 for convenience we refer to as "incoming"
 and "outgoing": each class consists of ends that are asymptotic
 to one another, and the decision of naming one class incoming, the other outgoing,
 is arbitrary.  We write $e\in {\sf in}(v)$ or $e\in {\sf out}(v)$ if $e$ has
 an end belonging to the appropriate class.  See Figure 2.  (Note: it can happen that
 $e$ belongs to both ${\sf in}(v)$ and ${\sf out}(v)$.)

  \begin{figure}[htbp]
\centering \epsfig{file=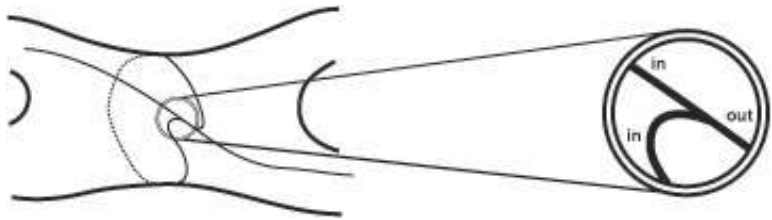, width=4in} \caption{Incoming and outgoing
ends. }
\end{figure}

 A {\em switch-additive measure} on $\tau$
 is a function $m:E\rightarrow\R_{+}$ for which
 \[ \sum_{e\in {\sf in}(v)} m(e) \;\; =\;\;\sum_{e\in {\sf out}(v)} m(e) \]
 for all switches $v$.
 The set of all switch-additive measures forms a linear cone ${\sf C}_{\tau}$
 in $\R^{E}$.
 
 Let $N(\tau )$ be a tubular neighborhood of $\tau$ equipped with a 
 (singular) foliation by line segments transverse to $\tau$.  See Figure 3.
 
  \begin{figure}[htbp]
\centering \epsfig{file=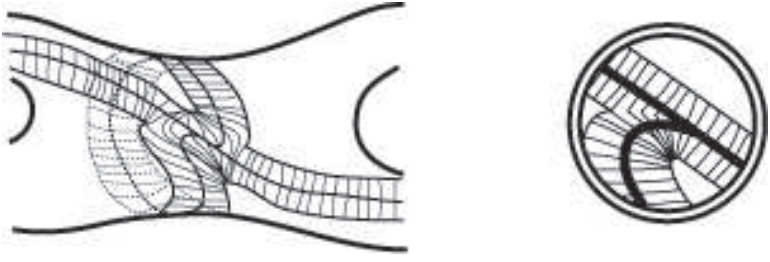, width=4in} \caption{The tubular
neighborhood $N(\tau )$. }
\end{figure}
 
A measured lamination $\mathfrak{f}\subset \Sigma$ is said to
be {\em carried} by $\tau$ if it may by isotoped into $N (\tau )$ transverse
to its foliation. We write in this case $\mathfrak{f}<\tau$. The subspace of 
isotopy classes of measured laminations
carried by $\tau$ is denoted $\mathcal{ML}_{\tau}(\Sigma )$.
There is an open surjection
\[ {\sf C}_{\tau}\longrightarrow \mathcal{ML}_{\tau}(\Sigma ) \]
which is a homeomorphism if $\tau$ is a train track.

Let $\Phi:\Sigma\rightarrow\Sigma$ be a pseudo Anosov diffeomorphism.  We say
that $\Phi$ {\em acts} on $\tau$ if $\Phi (\tau )$ may be isotoped into $N(\tau )$
transverse to its foliation.  We write then $\Phi (\tau )<\tau$.  Fix
a leaf $t_{i}\subset N(\tau )$ through each edge $e_{i}$ of $\tau$.  The
carrying matrix of $\Phi$ is by definition $M_{\Phi}=(a_{ij})$ where
\[ a_{ij}\;\; =\;\; |\Phi (e_{i})\cap t_{j}| . \]
$M_{\Phi}$ induces an inclusion 
${\sf C}_{\Phi (\tau )}\hookrightarrow {\sf C}_{\tau}$ which when precomposed
with the pushforward map ${\sf C}_{\tau}\rightarrow {\sf C}_{\Phi (\tau )}$
defines a linear map
\[ M_{\Phi}:{\sf C}_{\tau}\longrightarrow {\sf C}_{\tau}. \]

Note that the carrying matrix $M_{\Phi}$ is non-negative.  Such a matrix
has a unique eigenvalue of greatest modulus, which is positive-real and simple
\cite{Ga}.
This eigenvalue is called the {\em Perron root}.  A corresponding eigenvector
may be taken non-negative, and is called a {\em Perron vector}.  For $M_{\Phi}$,
the Perron root coincides with the entropy $\lambda$ of $\Phi$, and
the Perron vector parametrizes in track coordinates the unstable measured
lamination $\mathfrak{f}^{u}$ of $\Phi$.

When $\Phi$ is a pseudo Anosov of Thurston-Penner type, one can recover
the carrying matrix and all of its Perron data from a simpler
matrix which records the action of $\Phi$ on the curves in the families
$\mathcal{C}$ and $\mathcal{D}$.  Indeed, let $\tau$ be the bigon track
formed from $\mathcal{C}\cup\mathcal{D}$.  If $e_{i},e_{j}\in E_{\tau}$
 are edges contained in say $c$, $d$ resp.\ then given $e_{i'}\subset c$
 another edge, there exists a unique $e_{j'}\subset d$ with $a_{i'j'}=a_{ij}$.  
 Conversely, for $e_{j'}\subset d$, there exists $e_{i'}\subset c$ with
$a_{i'j'}=a_{ij}$.  It follows that the carrying matrix of $\Phi$ can be subdivided into
blocks indexed by pairs $(c,d)$, which are of the form
$a_{c,d}I$ where $I$ is a square matrix in which each column and row has
exactly one non zero entry = 1.  The matrix $(a_{c,d})$ whose columns and
rows are indexed
by $\mathcal{C}\cup \mathcal{D}$ has exactly the same Perron root as $M_{\Phi}$,
and its Perron vector gives that of $M_{\Phi}$ in the obvious way.  We
shall call this matrix the {\em curve matrix} of the Thurston-Penner
type pseudo Anosov $\Phi$, and we shall denote it $M_{\Phi}$ as well.

For example, if $\mathcal{C}=\{ c\}$, $\mathcal{D}=\{ d\}$ with $|c\cap d|=r$
and $\Phi=G^{-N}_{d}\circ
F^{N}_{c}$, then 
\[ M_{\Phi}\;\; =\;\;\left(
                      \begin{array}{ll}
                      1 & rN \\
                      rN & (rN)^{2}+1  
                      \end{array} \right) . \] 
We note that using the quadratic formula, it is easy to see that the eigenvalue
{\em not} equal to the Perron root is $<1$.

We now discuss tracks on the solenoid $\widehat{\Sigma}$: in fact,
we will only require tracks pulled back from surfaces appearing
in its defining inverse system.
Thus, if $\tau$ is a train track on such a surface $Z$, 
its preimage $\hat{\tau}$ is a smooth 1-dimensional branched solenoid 
with edge set $\widehat{E}\approx E\times \widehat{T}_{Z}$,
where $E$ is the edge set of $\tau$ and $\widehat{T}_{Z}$ is the fiber over a point of $Z$, homeomorphic
to the Cantor group $\hat \pi_{1}Z$.  With respect to this decomposition
we define a measure on $\widehat{E}$
by the formula
\[ \mu_{\widehat{E}}\;\; =\;\; \mu_{E}\times \eta_{\widehat{T}_{Z}} \]
where $\mu_{E}(e)=1$ for each edge of $E$ and $\eta_{\widehat{T}_{Z}}$
is the restriction to $\widehat{T}_{Z}$ of the transverse invariant measure of $\widehat{\Sigma}$.
In addition, if $\tau$ is equipped with a 
switch additive measure $\upsilon$, the pullback $\hat{\upsilon}$
is a transversally continuous switch-additive measure
on $\hat{\tau}$; the cone of such measures on $\hat{\tau}$ is denoted ${\sf C}_{\hat{\tau}}$.
The relation $\hat{\mathfrak{f}}<\hat{\tau}$
has exactly the same meaning as in the case of a surface.

\section{Intersection Formulas}

Let $\tau ,\kappa$ be bigon tracks in $\Sigma$ that intersect transversally
and minimally with edge sets $E_{\tau}$ and $E_{\kappa}$; 
let $\mathfrak{f}, \mathfrak{g}$ be measured laminations
carried by them, parametrized by weights $\upsilon ,\omega $.  The intersection
pairing may be calculated by the following formula \cite{Pe1}:
\begin{equation}\label{intersection} \text{\bf I}( \mathfrak{f}, \mathfrak{g} )
\;\; =\;\; 
\sum_{ e\in E_{\tau},\;   e'\in E_{\kappa }}
\upsilon (e)\omega(e') .\end{equation}
It is useful to re-express (\ref{intersection}) as a sum over edges
in $E_{\tau} $ only.  Thus if we write
\[\omega (e)\;\; =\;\;\sum_{e'\in E_{\kappa }}\omega (e')|e\cap e'|\]
then 
\begin{equation}\label{intersection2}
 \text{\bf I}( \mathfrak{f}, \mathfrak{g} )\;\; =\;\;
 \sum_{e\in E_{\tau }}\upsilon (e)\omega (e).
 \end{equation}

 Suppose now that $\hat{\mathfrak{f}}$, 
 $\hat{\mathfrak{g}}$ are measured laminations obtained as preimages of
 measured laminations $\mathfrak{f}\subset Y$ and $\mathfrak{g}\subset Z$,
 surfaces occurring in the defining system of $\widehat{\Sigma}$.  Let
 $W$ be a surface finitely covering each of $Y,Z$, and let $\tilde{\mathfrak{f}}$,
$\tilde{\mathfrak{g}}$ be the preimages in $W$
 of $\mathfrak{f}$, $\mathfrak{g}$.  Let $\deg (W)$ be the degree of
 the covering $W\rightarrow \Sigma$.

\begin{prop}  
$\text{\bf I}( \hat{\mathfrak{f}}, 
 \hat{\mathfrak{g}} )\;\; =\;\; 
 \frac{\text{\bf I}(\tilde{\mathfrak{f}}, \tilde{\mathfrak{g}})}{\deg (W)}$.
\end{prop} 

\begin{proof} The intersection locus of $\hat{\mathfrak{f}}$ and 
 $\hat{\mathfrak{g}}$ is of the form
 \[ (\tilde{\mathfrak{f}}\cap \tilde{\mathfrak{g}})\times \widehat{T}_{W} \]
where $\widehat{T}_{W}$ is a fiber of $\widehat{\Sigma}\rightarrow W$.  Since
$\widehat{T}_{W}$ has $\eta$-measure $1/\deg (W)$, the result follows.
\end{proof}
 
Let $\hat{\mathfrak{f}}, 
 \hat{\mathfrak{g}}$ be as in the previous paragraphs.  Suppose now that $\mathfrak{f}$
 is parametrized by $\upsilon :E\rightarrow\R_{+}$ a weight on a bigon track
 $\tau\subset Y$ and $\mathfrak{g}$ is parametrized by 
 $\omega :E'\rightarrow\R_{+}$ a weight on a bigon track
 $\tau'\subset Z$.  The preimages $\tilde{\mathfrak{f}}, \tilde{\mathfrak{g}}$
 are parametrized by the pullback weights $\tilde{\upsilon}$, $\tilde{\omega}$ on the preimages
 $\tilde{\tau},\tilde{\tau}'$.  Rewriting $\omega$ as above
as a function of the edge set $E$, we have
\[ \text{\bf I}(\tilde{\mathfrak{f}}, \tilde{\mathfrak{g}})\;\; =\;\;
\sum_{\tilde{e}\in\widetilde{E}}\tilde{\upsilon}(\tilde{e})\tilde{\omega}(\tilde{e}) . \]
Let $\hat{\tau}$ be the preimage of $\tau$ in $\widehat{\Sigma}$, and let
$\hat{\upsilon}$, $\hat{\omega}$ be the pullbacks along the projection $\widehat{E}\rightarrow
\widetilde{E}$.

\begin{prop}  $\text{\bf I}( \hat{\mathfrak{f}}, 
 \hat{\mathfrak{g}} )\;\; =\;\; 
 \int_{\widehat{E}}\hat{\upsilon}\hat{\omega}\; d\mu_{\widehat{E}}$ where
 $\mu_{\widehat{E}}$ is the edge measure on $\widehat{E}$.
\end{prop}

\begin{proof}  A calculation:
\[
\text{\bf I}( \hat{\mathfrak{f}}, \hat{\mathfrak{g}} ) 
\;\; = \;\;
\left( \sum_{\tilde{e}\in\widetilde{E}}
\tilde{\upsilon}(\tilde{e})\tilde{\omega}(\tilde{e})\right)\cdot\frac{1}{\deg (W)} \\
\;\; = \;\;
\left(\sum_{\tilde{e}\in\widetilde{E}}
\tilde{\upsilon}(\tilde{e})\tilde{\omega}(\tilde{e})\right) |\widehat{T}_{W}| \\
\;\; = \;\;
 \int_{\widehat{E}}\hat{\upsilon}\hat{\omega}\; d\mu_{\widehat{E}}.
\]
\end{proof}

\section{Teichm\"{u}ller Theory of $\widehat{\Sigma}$}

References for material in this section are \cite{Ge1}, \cite{Na-Su}, \cite{Su}.

The definition of the Teichm\"{u}ller space $\mathcal{T}(\widehat{\Sigma})$
and of its metric $d_{\mathcal{T}(\widehat{\Sigma})}$
copies that of a surface.  In particular, 
\begin{enumerate}
\item A conformal structure on $\widehat{\Sigma}$ is determined by a conformal
structure on each leaf.  These structures are required to
vary continuously in the transverse direction.
\item Elements of $\mathcal{T}(\widehat{\Sigma})$ are represented 
by marked solenoids {\em i.e.}\ by homeomorphisms
$\mu :\widehat{\Sigma}\rightarrow\widehat{\Sigma}_{\mu}$, where $\widehat{\Sigma}_{\mu}$
is presumed to have a conformal structure. 
\item The marked solenoid
$\mu' :\widehat{\Sigma}\rightarrow\widehat{\Sigma}_{\mu'}$ is equivalent
to $\mu$ if there exists an isomorphism $\sigma :\widehat{\Sigma}_{\mu}\rightarrow
\widehat{\Sigma}_{\mu'}$ such that $\sigma\circ\mu\simeq\mu'$. 
\end{enumerate}

$\mathcal{T}(\widehat{\Sigma})$ has the structure of a separable Banach manifold.
The canonical projection
$\hat{p}:\widehat{\Sigma}\rightarrow Z$ onto any surface $Z$ in the
defining inverse system induces a direct system of 
isometric inclusions $\hat{p}^{\ast}:\mathcal{T}(Z)\hookrightarrow
\mathcal{T}(\widehat{\Sigma})$.  

\begin{theo}[\cite{Na-Su}]\label{teichinfty}  The induced inclusion
\[ i:\lim_{\longrightarrow}\mathcal{T}(Z)\;\hookrightarrow\; 
\mathcal{T}(\widehat{\Sigma})\]
is isometric with dense image.
\end{theo}

For most of our purposes, it will be sufficient to work with the dense subspace 
\[ \mathcal{T}^{\infty}\;\; :=\;\; i(\lim_{\longrightarrow}\mathcal{T}(Z)) ,\]
which is an incomplete metric space with respect to the direct limit of
the Teichm\"{u}ller metrics and a pre-Banach manifold.
Unless otherwise said, all structures $\hat{\mu}$ considered below will be
assumed to be in $\mathcal{T}^{\infty}$.  

Let $\widehat{\Sigma}_{\hat{\mu}}$ be as above.
By a holomorphic quadratic differential $\hat{q}$ on $\widehat{\Sigma}_{\hat{\mu}}$, 
we shall always mean the pull-back of a holomorphic quadratic differential $q$ occurring
on some surface $Z_{\mu}$, where $\hat{\mu}$ is the pull-back of $\mu$.  
Thus, $\hat{q}$ is 
a choice of holomorphic quadratic differential on each leaf, constant along
the fiber transversals $\widehat{T}_{Z}$ over $Z$.
The tangent space to 
$\mathcal{T}^{\infty}$ at $\hat{\mu}$ may be identified with the direct limit
\[  Q_{\hat{\mu}}^{\infty}\;\; =\;\;\lim_{\longrightarrow}Q_{\mu}(Z) .\]
The
tangent bundle of $\mathcal{T}^{\infty}$ is 
then identified with $\mathcal{Q}^{\infty}=\lim_{\rightarrow} \mathcal{Q}(Z)$
where $\mathcal{Q}(Z)$ is the space of holomorphic quadratic differentials
on $Z$ (with respect to all possible complex structures).

The $L^{1}$ norm on $Q^{\infty}_{\hat{\mu}}$ is defined
\[  \| \hat{q}\|\;\; =\;\; \int_{\widehat{\Sigma}_{\mu}}|\hat{q}|\, d\eta .\]
From
the pre-Finsler norm $\|\cdot\|$ we induce 
a path metric $d_{L^{1}}$ on $\mathcal{T}^{\infty}$ that
defines 1) the $L^{1}$ topology on $\mathcal{T}^{\infty}$ and 2)
along with $\|\cdot\|$, the $L^{1}$ topology 
on $\mathcal{Q}^{\infty}$.  

To any quadratic diferential $\hat{q}$ one associates two transverse, measured laminations
$\hat{\mathfrak{f}}^{h}$ and $\hat{\mathfrak{f}}^{v}$ (those that correspond
to the horizontal and vertical
trajectories of $\hat{q}$).  We have the following generalization to $\widehat{\Sigma}$
of a well-known formula for surfaces:
\begin{lemm}$ \text{\bf I} (\hat{\mathfrak{f}}^{h},\hat{\mathfrak{f}}^{v}) =\|\hat{q}\|$ .
\end{lemm}

\begin{proof} $\hat{q}$ is the lift of a holomorphic quadratic differential
on some surface $Z_{\mu}$, where $\hat{\mu}$ is the lift of $\mu$.  The result
now follows from the classical formula and the fact that the lift of an area
measure on $Z_{\mu}$ scales by $1/(\deg Z )$ in $\widehat{\Sigma}$.
\end{proof}

In the same way,
we may also avail ourselves of a direct limit version of the theorem of Hubbard
and Masur \cite{HuMa}:

\begin{theo}\label{HubbMas} Any pair of measured laminations 
$\hat{\mathfrak{f}},\hat{\mathfrak{g}}\in\mathcal{ML}^{\infty}$
determines a unique quadratic differential $\hat{q}$.
\end{theo}

\begin{theo}\label{intconvimpl1}  Let $\hat{q},\hat{q}_{i}\in\mathcal{Q}^{\infty}$, $i=1,2,\dots$, 
be quadratic differentials.
If $\hat{\mathfrak{f}}_{i}^{h}\rightarrow\hat{\mathfrak{f}}^{h}$ and
$\hat{\mathfrak{f}}_{i}^{v}\rightarrow\hat{\mathfrak{f}}^{v}$
in the intersection topology,
then $\hat{q}_{i}\rightarrow \hat{q}$ in the $L^{1}$ topology.
\end{theo}

\begin{proof} We assume first that there exists 
$\hat{\mu}\in\mathcal{T}^{\infty}$ with
$\hat{q}_{i},\hat{q}\in Q^{\infty}_{\hat{\mu}}$.  
Let $\hat{\sf f}^{h}$, $\hat{\sf f}^{v}$ and 
$\hat{\sf f}^{h}_{i}$, $\hat{\sf f}_{i}^{v}$ be the pairs of measured foliations
which are the horizontal and vertical line fields of $\hat{q}$ resp.\ $\hat{q}_{i}$.  
By the comments of \S 3, the hypothesis on the convergence of measured laminations
is equivalent to the corresponding statement for measured foliations.
In particular we have, 
\begin{equation}\label{intconv1}
\lim \text{\bf I}(\hat{\sf f}^{h},\hat{\sf f}^{h}_{i})\;\;=\;\; 0
\quad\quad\text{and}\quad\quad
\lim \text{\bf I}(\hat{\sf f}^{v},\hat{\sf f}^{v}_{i})\;\;=\;\; 0 .
\end{equation}
Consider smooth measured foliations $\hat{\sf g}_{i}^{h}$,
$\hat{\sf g}_{i}^{v}$ equivalent to $\hat{\sf f}_{i}^{h}$,
$\hat{\sf f}_{i}^{v}$ whose heights with respect to $\hat{\sf f}^{h}$,
$\hat{\sf f}^{v}$ nearly give the intersections 
$\text{\bf I}(\hat{\sf f}^{h},\hat{\sf f}^{h}_{i})$,
$\text{\bf I}(\hat{\sf f}^{v},\hat{\sf f}^{v}_{i})$.  More precisely,
for $\epsilon_{i}\rightarrow 0$,
\begin{equation}\label{intconv2}
\int_{\hat{\sf g}_{i}^{h}} |{\sf Im}\sqrt{\hat{q}}|\; d\eta  
-  \text{\bf I}(\hat{\sf f}^{h},\hat{\sf f}^{h}_{i}) \;\; <\;\; \epsilon_{i}
\quad\quad \text{ and } \quad\quad
\int_{\hat{\sf g}_{i}^{v}} |{\sf Re}\sqrt{\hat{q}}|\; d\eta  
-  \text{\bf I}(\hat{\sf f}^{v},\hat{\sf f}^{v}_{i}) \;\; <\;\; \epsilon_{i} .
\end{equation}
Let $\hat{q}_{i}'$ denote the smooth quatratic differential
whose
horizontal and vertical foliations 
are $\hat{\sf g}^{h}_{i}$, $\hat{\sf g}^{v}_{i}$. 
Given $\delta >0$,
let $\hat{A}_{i}$ be the set of points for which 
$|\hat{q}-\hat{q}_{i}'|$ is uniformly $\delta$ small.
Let $\hat{B}_{i}=
\widehat{\Sigma}\setminus \hat{A}_{i}$.

We begin by showing that  
\begin{equation}\label{zeroconv}
\int_{\hat{B}_{i}}|\hat{q}-\hat{q}'_{i}|\; d\eta \rightarrow 0
\end{equation}
as $i\rightarrow 0$.
By (\ref{intconv1}) and (\ref{intconv2}),
$\int_{\hat{B}_{i}}|\hat{q}|\,d\eta \rightarrow 0$.
If for $i$ large, there is a $m>0$ with
\[ 0\;\;<\;\; m\;\;\leq\;\; \int_{\hat{B}_{i}}|\hat{q}-\hat{q}_{i}'|\; d\eta\]
we must also have
\begin{equation}\label{posqimass}
0\;\;<\;\;  m_{0}\;\;\leq\;\; 
 \int_{\hat{B}_{i}}|\hat{q}_{i}'|\; d\eta 
 \end{equation}
 for some $m_{0}>0$.
The fact that $\int_{\hat{B}_{i}}|\hat{q}|\, d\eta\rightarrow 0$ whereas
$\int_{\hat{B}_{i}}|\hat{q}_{i}'|\, d\eta$ does not would violate
(\ref{intconv1}) and (\ref{intconv2}) as well.
Thus 
$\int_{\hat{B}_{i}}|\hat{q}_{i}|\; d\eta \rightarrow 0$,
proving (\ref{zeroconv}).
In particular, we have
\[  \lim \int |\hat{q}-\hat{q}'_{i}|\; d\eta \;\; =\;\; 
\lim \int_{ \hat{A}_{i} } |\hat{q}-\hat{q}'_{i}| .\]

Now $\hat{q}_{i}'$ is measure
equivalent to $\hat{q}_{i}$, hence 
$\hat{q}-\hat{q}_{i}$ is measure equivalent to $\hat{q}-\hat{q}'_{i}$. 
By the second minimal norm property \cite{Gar}, it follows that
\[  \|\hat{q}-\hat{q}_{i}\|\;\leq\; \|\hat{q}-\hat{q}_{i}'\| . \]
Letting $\epsilon\rightarrow 0$, 
we obtain $\|\hat{q}-\hat{q}_{i}\|\rightarrow 0$.
This proves the theorem in the special case where $\hat{\mu}=\hat{\mu}_{i}$
for $i$ large.

Now we suppose that $\hat{q}\in Q^{\infty}_{\hat{\mu}}$, 
$\hat{q}_{i}\in Q^{\infty}_{\hat{\mu}_{i}}$
with $\hat{\mu}_{i}\not= \hat{\mu}$.
Let $\hat{C}_{i}$ be the set of points where 
the foliations 
$\hat{\sf f}_{i}^{h}$,
$\hat{\sf f}_{i}^{v}$ are uniformly $\epsilon$-close to $\hat{\sf f}^{h}$,
$\hat{\sf f}^{v}$ 
in the $\hat{q}$-metric.  Then
for $i$ large, in $\hat{C}_{i}$
the complex structures
defined by $\hat{q}_{i},\hat{q}$ are uniformally nearly conformal.  On the other
hand, in $\hat{D}_{i}=\widehat{\Sigma}\setminus \hat{C}_{i}$ they are not, 
but the $|\hat{q}|$-volume of this set 
limits to zero.  Therefore, if $\hat{p}_{i}\in  Q^{\infty}_{\hat{\mu}}$ 
 generates the Teichm\"{u}ller geodesic connecting
$\hat{\mu}$ to $\hat{\mu}_{i}$ in time 1, it follows that 
$\| \hat{p}_{i}\|\rightarrow 0$
so that $\hat{\mu}_{i}$ converges to $\hat{\mu}$ in the $L^{1}$ path metric. 
Moreover, the induced flow of quadratic differentials takes 
$\hat{q}$ $L^{1}$ close to $\hat{q}_{i}$ so that
$\hat{q}_{i}$ converges to $\hat{q}$ in the $L^{1}$ topology
on $\mathcal{Q}^{\infty}$. 
\end{proof}

\section{The $L^{1}$ Ehrenpreis Conjecture}

The $L^{1}$ {\sf EC} is the following statement:

\begin{l1ec}  The mapping class group ${\sf M}(\widehat{\Sigma}, \ell )$ 
acts with $L^{1}$ dense
orbits on $\mathcal{T}(\widehat{\Sigma})$.
\end{l1ec}

Since $\mathcal{T}^{\infty}$ is dense in $\mathcal{T}(\widehat{\Sigma})$
and ${\sf M}(\widehat{\Sigma}, \ell )$ stabilizes $\mathcal{T}^{\infty}$ \cite{Ge1},
it will be enough to demonstrate that ${\sf M}(\widehat{\Sigma}, \ell )$ acts with $L^{1}$ dense
orbits on $\mathcal{T}^{\infty}$.

In \cite{Ma}
it is shown that for a closed surface $Z$, a Teichm\"{u}ller dense subset of 
pairs $\mu , \nu\in \mathcal{T}(Z)$ lie
on the axes of pseudo Anosov homeomorphisms.  
By its definition as an isometric direct limit, $\mathcal{T}^{\infty}$
enjoys the same property.  Fix a pair $\hat{\mu},\hat{\nu}\in \mathcal{T}^{\infty}$;
without loss of generality, we may then assume that $\hat{\mu},\hat{\nu}$
lie on the axis $A$ of a pseudo Anosov diffeomorphism $\hat{\Phi}$ which is the
lift of a pseudo Anosov $\Phi:Z\rightarrow Z$, for some surface $Z$ occurring
in the defining system of $\widehat{\Sigma}$. 

By an {\em $L^{1}$ ${\rm n}$th root} of $\hat{\Phi}$ we mean a sequence $\{ \hat{\Psi}_{m}\}$
of pseudo Anosov homeomorphisms for which
\begin{enumerate}  
\item If $\lambda_{m}$ is the entropy of $\hat{\Psi}_{m}$ then
$\lim\lambda_{m}=\lambda^{1/n}$.
\item  If $A_{m}$ is the axis of $\hat{\Psi}_{m}$ then $A_{m}\rightarrow A$
converges in the Hausdorff topology induced from the $L^{1}$ metric.
\end{enumerate}

\begin{theo}  If for every pseudo Anosov $\hat{\Phi}$ and every $n>0$,
$\hat{\Phi}$ has an $L^{1}$ $n$th root, 
then the $L^{1}$ {\sf EC} is true.
\end{theo}

\begin{proof}  Suppose that $\hat{\mu},\hat{\nu}$ lie on the
axis $A$ of $\hat{\Phi}$ and
let $\{\hat{\Psi}_{m}\}$ be an $L^{1}$ $n$th root, $n$ large.  
Since the $A_{m}\rightarrow A$ in the $L^{1}$
Hausdorff topology, there exists $\hat{\mu}', \hat{\nu}'$ lying on some axis $A_{m}$
with $d_{L^{1}}(\hat{\mu} ,\hat{\mu}' )<\epsilon $, 
$d_{L^{1}}(\hat{\nu} ,\hat{\nu}' )<\epsilon$.  On the
axis $A_{m}$, we may move via a power of $\hat{\Psi}_{m}$ $\hat{\mu}'$ close to 
$\hat{\nu'}$, which
implies that $\hat{\mu}$ is moved close to $\hat{\nu}$ by the same power of
$\hat{\Psi}_{m}$ as well.
\end{proof}

\begin{note}  The existence of {\em Teichm\"{u}ller roots} for all $\hat{\Phi}$
and all $n$ implies
the classical ${\sf EC}$.
\end{note}

\section{Directional Density}

A family $\mathcal{P}$ of pseudo Anosov homeomorphisms is said to
be {\em directionally dense} (in $\mathcal{Q}^{\infty}$) if the set of quadratic differentials
tangent to axes of elements of $\mathcal{P}$ is
Teichm\"{u}ller dense in $\mathcal{Q}^{\infty}$.
By \cite{Ma}, the family of all pseudo Anosov maps
is directionally dense.  In fact, 
it follows easily from
the arguments in \cite{Ma} that the family of lifts of pseudo Anosovs $\Phi$
of the type $\Phi =G^{-2N}\circ F^{2N}$, where
$F$, $G$
are right Dehn twists about simple closed curves $c, d$ that 
fill $Z$, where $Z$ ranges over all surfaces in the 
defining system of $\widehat{\Sigma}$,
is directionally
dense.  By Corollary 2.6 in \cite{Ma}, the subfamily obtained
by demanding that $c$
is nonseparating is also directionally dense.

Now given a nonseparating simple closed curve $\gamma\subset Z$,
let $\rho_{\gamma}:Z_{\gamma}\rightarrow Z$ be the degree 2 cover
obtained by cutting two copies of $Z$ along $\gamma$ and gluing 
ends.  We say that a pair $(c,d)$ of filling, simple closed curves
is {\em interlacing} if there exists a pair of
nonseparating simple closed curves $\alpha ,\beta$ such that 
$\rho^{-1}_{\alpha}(c)$ is connected whereas $\rho^{-1}_{\alpha}(d)$ is not
and $\rho^{-1}_{\beta}(d)$ is connected whereas $\rho^{-1}_{\beta}(c)$
is not.  See Figure 4.

 \begin{figure}[htbp]
\centering \epsfig{file=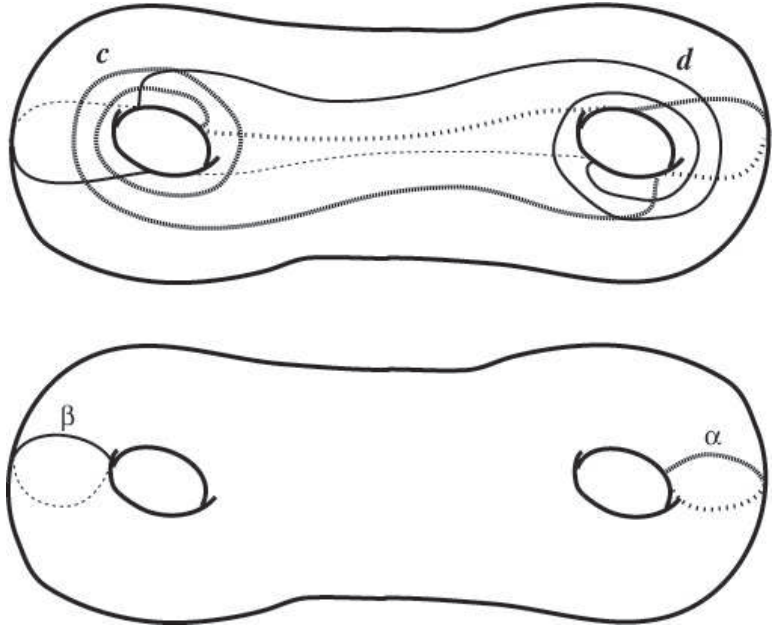, width=4in} \caption{A filling pair $c,d$
and a pair $\alpha ,\beta$ interlacing them.}
\end{figure}

 Let $\mathcal{P}$ be the family 
of pseudo Anosov homeomorphisms
of $\widehat{\Sigma}$ which
are lifts of pseudo Anosovs of the form
\[ \Phi = G^{-2N}\circ F^{2N}:Z\rightarrow Z,\] 
where 
\begin{enumerate}
\item $F$, $G$
are right Dehn twists about $c,d\subset Z$.
\item $c$ is nonseparating and 
$(c,d)$ is an interlacing pair. 
\item $Z$ ranges over all surfaces in
the defining system of $\widehat{\Sigma}$.
\end{enumerate}

\begin{lemm}\label{lacing} $\mathcal{P}$
is directionally
dense.
\end{lemm}

\begin{proof}  It is enough to show that for a fixed surface $Z$,
the family of maps satisfying (1) and (2) is directionally dense in $\mathcal{Q}(Z)$.  
Assume that $c$ and $d$ are filling, generating the pseudo
Anosov homeomorphism $\Phi= G^{-2N}\circ F^{2N}$.  If $(c,d)$ is not an interlacing pair,
there exists a simple closed curve
$\delta$ for which the pair $(c,\delta )$ is interlacing, though
not necessarily filling.  Indeed, one may assume after a homeomorphism
that $c$ is the curve appearing in Figure 5; then taking $\delta,\alpha, \beta$ as 
indicated there, $(c,\delta )$ is interlacing with respect to the pair
$(\alpha$, $\beta )$.   Now for $j$ large, $\delta_{j}=G^{j}(\delta )$
is close to $d$, hence $(c,\delta_{j})$ is eventually filling.  
If $j$ is in addition
even, $G^{j}$ lifts to the total space of any degree 2 cover of $\Sigma$,
thus the pair $(c,\delta_{j})$ is interlacable with respect to the same
curves interlacing $(c,\delta )$.  For $j$ large, the pseudo Anosov
$\Phi_{j}=G^{-2N}_{\delta_{j}}\circ F^{2N}$ has axis close to that of $\Phi$, 
and since the maps of the form
$\Phi$ are already directionally dense, we are done.
\end{proof}

 \begin{figure}[htbp]
\centering \epsfig{file=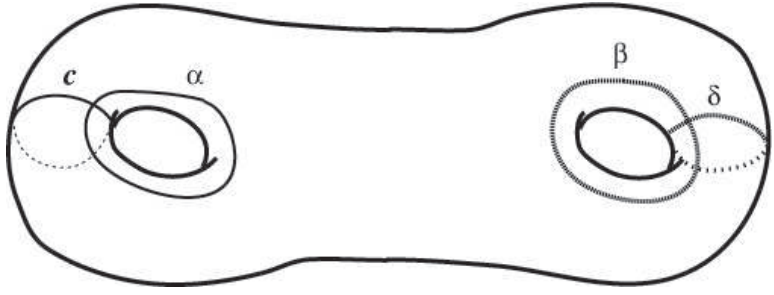, width=4in} \caption{Every nonseparating
$c$ is a member of a (not necessarily filling) interlacing pair.}
\end{figure}

\section{Necklace Roots}

Let $n\in\N$ and let $\hat{\Phi}\in\mathcal{P}$, the family appearing in 
Lemma~\ref{lacing}, so that 
in particular $\hat{\Phi}$ is the lift of a pseudo Anosov 
of the form $\Phi=G^{-2N}\circ F^{2N}: Z\rightarrow Z$.  
Let $\rho_{mn}:Z_{mn}\rightarrow Z$
be the cover obtained by cutting $2mn$ copies of $Z$ 
along a pair $\alpha$ and $\beta$ interlacing $c,d$ and gluing in a circular
fashion.  We call $\rho_{mn}$ the {\em necklace cover} associated to $(c,d)$. 
In Figure 6, we illustrate the construction of the necklace $Z_{mn}$ 
and the formation of the
lifts of the curve $c$.  In Figure 7 we display the finished necklace.

 \begin{figure}[htbp]
\centering \epsfig{file=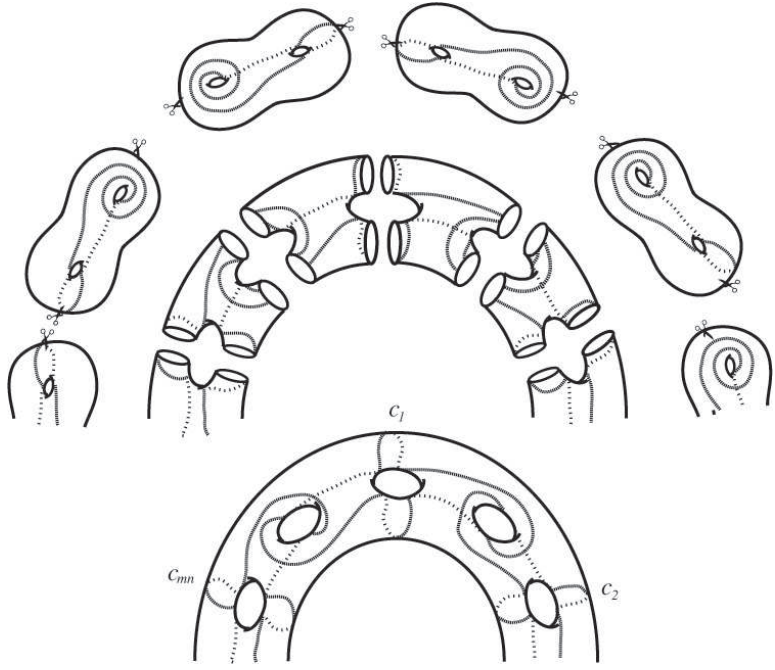, width=4in} \caption{Making the necklace. }
\end{figure}

\begin{figure}[htbp]
\centering \epsfig{file=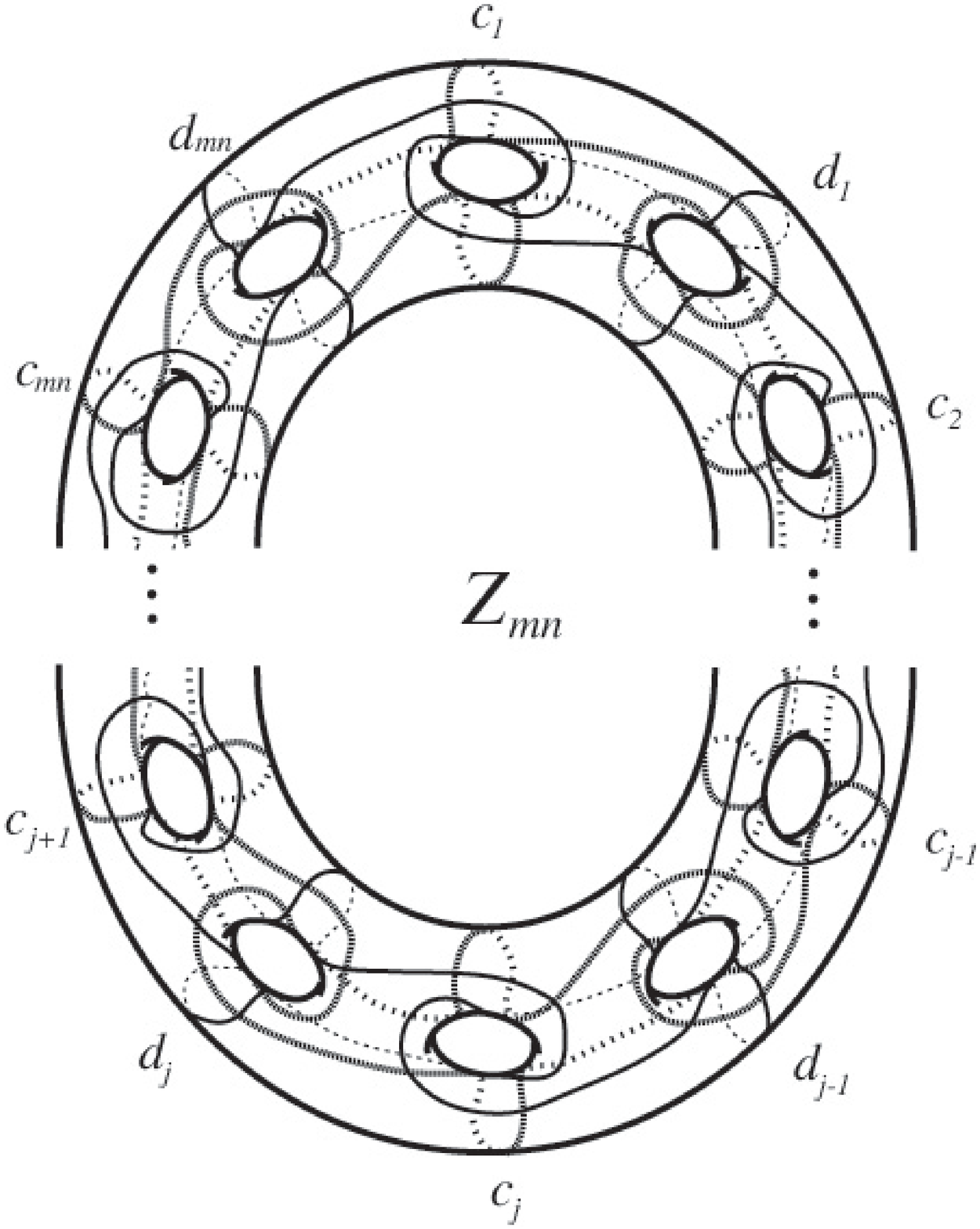, width=4in} \caption{The finished necklace. }
\end{figure}

There are $mn$ lifts $c_{1},\dots ,c_{mn}$
and $d_{1},\dots ,d_{mn}$ of each of $c$ and $d$, each mapping with
degree two onto their ancestor.  On $Z_{mn}$, 
$\Phi$ lifts to 
\[\tilde{\Phi}\;\; =\;\; G_{mn}^{-N}\circ\dots\circ G_{1}^{-N}\circ
F_{mn}^{N}\circ\dots\circ F_{1}^{N}   \]
where $F_{i}$, $G_{i}$ is the right Dehn twist about $c_{i}$, $d_{i}$.
Let $\chi$ denote the clockwise rotation of $Z_{mn}$ by an angle
of $2\pi /n$, so that the pair $c_{i}$, $d_{i}$ is taken to
$c_{j+m}$, $d_{j+m}$ (indices taken mod $mn$).  We define
the {\em necklace ${\rm n}$th root} to be the sequence
of lifts of pseudo Anosovs
$\{\widehat{\sqrt[n]{\Phi}_{m}}\}$ to $\widehat{\Sigma}$ where
\[ \sqrt[n]{\Phi}_{m}\;\; =\;\; \chi\circ  G_{m}^{-N}\circ\dots\circ G_{1}^{-N}\circ
F_{m}^{N}\circ\dots\circ F_{1}^{N} ,\]
$m=2,3,\dots $.  The necklace $n$th root is the basic contruction used 
in the formation of $L^{1}$ roots.  The construction of
$\sqrt[n]{\Phi}_{m}$ is a generalization of one
that first appeared in \cite{Pe3}, where branched covers
were used.

\begin{lemm}  $\sqrt[n]{\Phi}_{m}$ is pseudo Anosov for all $m$.
\end{lemm}
 
 \begin{proof}  For $i=1,\dots ,n$, let
 \[ T_{i}\;\; =\;\; G_{im}^{-N}\circ\dots\circ G_{(i-1)m+1}^{-N}\circ
 F_{im}^{N}\circ\dots\circ F^{N}_{(i-1)m+1}.  \]
 Then it is easy to see that
 \[ (\sqrt[n]{\Phi}_{m})^{n}\;\; =\;\; T_{2}\circ \dots \circ T_{n}\circ T_{1},\]
 which is of Thurston-Penner type, hence
\cite{Pe1} $(\sqrt[n]{\Phi}_{m})^{n}$ is pseudo Anosov, 
 implying $\sqrt[n]{\Phi}_{m}$ is
 pseudo Anosov as well.
 \end{proof}
 
 \section{Existence of $L^{1}$ Roots}
 
Denote by $\mathcal{C}_{Z}(\widehat{\Sigma})$ the
family of simple closed solenoids which are lifts of simple closed
curves on $Z$. 
We begin by
constructing a family 
$\{ \hat{\Psi}_{m}\}$ whose stable and unstable laminations 
intersection converge to those of $\hat{\Phi}$ with respect to
test solenoids in $\mathcal{C}_{Z}(\widehat{\Sigma})$.   
 
For each $m=2,3,\dots$ let $\sqrt[mn]{\Phi}_{m}$ denote
the $mth$ element in the sequence of pseudo Anosovs
whose lifts define the $mn$th necklace root of $\hat{\Phi}$.
Define the sequence $\{ \hat{\Psi}_{m}\}$ as 
the lifts to $\widehat{\Sigma}$ of the pseudo Anosov homeomorphisms
 \[  \Psi_{m}\;\; =\;\; \big( \sqrt[nm]{\Phi}_{m}\big)^{m}.\]
 Observe that the stable and unstable foliations
 of $\Psi_{m}$ and $\sqrt[nm]{\Phi}_{m}$ are equal.
 
 \begin{note}
 $\hat{\Psi}_{m}$ is not the same as $\widehat{\sqrt[n]{\Phi}_{m}}$.
 In fact, if we lift $\sqrt[n]{\Phi}_{m}$ to $Z_{m^{2}n}$
 -- where $\Psi_{m}$ is defined -- we see that
 this lift twists along $m$ disjoint ``blocks'' of curves, each block
 consisting of a succession of $m$ lifts of $c$ and $d$.  On the other
 hand, $\Psi_{m}$ consists of twists along one block consisting of a succession
 of $m^{2}$ lifts of $c$ and $d$. 
 As we shall see, the stable and unstable laminations
 of the family $\{ \hat{\Psi}_{m} \} $ have better convergence properties
 than those of $\{ \widehat{\sqrt[n]{\Phi}_{m}} \} $.
 \end{note} 
 
 Denote by $\hat{\mathfrak{f}}_{m}^{u}, \hat{\mathfrak{f}}_{m}^{s}$
 and by $\hat{\mathfrak{f}}^{u}, \hat{\mathfrak{f}}^{s}$ the unstable
 and stable laminations of $\hat{\Psi}_{m}$ and $\hat{\Phi}$.
 
\begin{lemm}\label{Zintconv}  For all $\hat{c}\in\mathcal{C}_{Z}(\widehat{\Sigma})$, 
\[ \text{\bf I}(\hat{\mathfrak{f}}_{m}^{u},\hat{c})\longrightarrow
\text{\bf I}(\hat{\mathfrak{f}}^{u},\hat{c}) \quad\quad\text{ and }\quad\quad
\text{\bf I}(\hat{\mathfrak{f}}_{m}^{s},\hat{c})\longrightarrow
\text{\bf I}(\hat{\mathfrak{f}}^{s},\hat{c}) .
\] 
\end{lemm}

\begin{proof}  The proof will be through examination of curve matrices.
We begin with $\Phi$.  Let $r=\text{\bf I}(c,d)$. The action of $\Phi$ along
the curves $c,d$ is given by the matrix
\[  M_{\Phi}\;\; =\;\; \left( \begin{array}{cc}
                               1 & 2rN \\
                               2rN & (2rN)^{2}+1
                               \end{array}
  \right) . \]
Let $Z_{m^{2}n}$ be the surface where 
$\Psi_{m}$ is defined.  The curve families
$\mathcal{C}=\{ c_{1},\dots ,c_{m^{2}n}\}$, $\mathcal{D}=\{d_{1},\dots d_{m^{2}n}\}$ 
are filling, and
the action of $\tilde{\Phi}$ on $\mathcal{C}\cup\mathcal{D}$ is prescribed schematically by
the matrix in Figure 8, where all entries not contained in the boxed vectors
are zero, and where the ``broken'' vectors indicate that only that
portion of the corresponding vector is used.  For example, in the upper right
hand corner we have the entry $rN$, which is the bottom half of the $B$-vector; 
in the lower right hand corner, we have the vector entry 
$\big((rN)^{2}\;\; 2(rN)^{2}+1\big)^{\sf T}$,
which is the top two thirds of the vector $D$, and so on.
The curve matrix of $(\Psi_{m})^{n}=(\sqrt[mn]\Phi_{m})^{mn}$ is displayed
in Figure 9. The black vectors indicate regions where
$M_{(\Psi_{m})^{n}}=(M_{\Psi_{m}})^{n}=(M_{\sqrt[mn]\Phi_{m}})^{mn}$ 
differs from $\widetilde{M}_{\Phi}$.
Figure 10 contains the curve matrix of $\sqrt[mn]\Phi_{m}$.

\begin{figure}[htbp]
\centering \epsfig{file=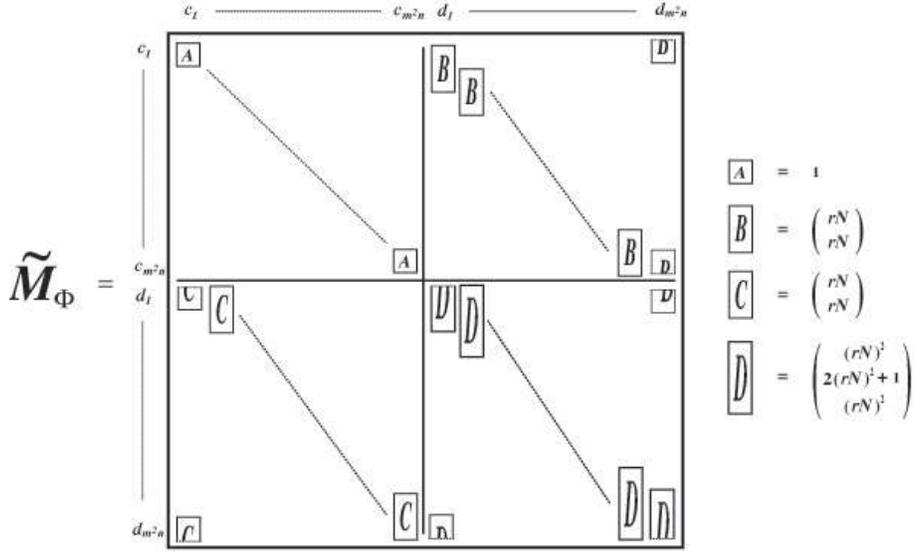, width=5.5in} \caption{The curve matrix
for $\tilde{\Phi}$. }
\end{figure}

\begin{figure}[htbp]
\centering \epsfig{file=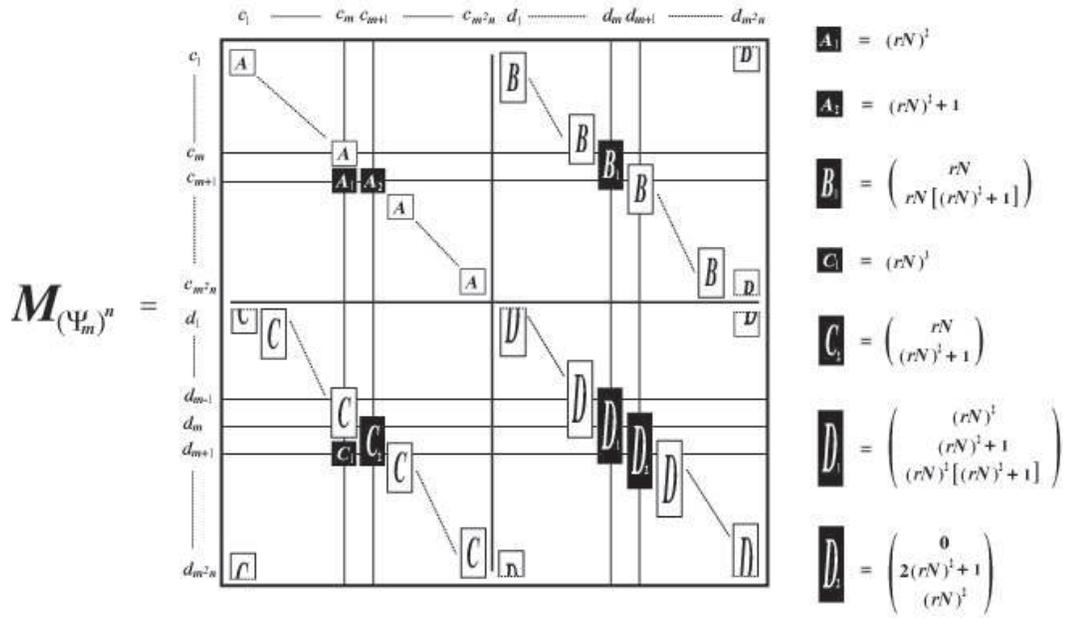, width=5.5in} \caption{The curve matrix
for $\Psi_{m}^{n}$. }
\end{figure} 

\begin{figure}[htbp]
\centering \epsfig{file=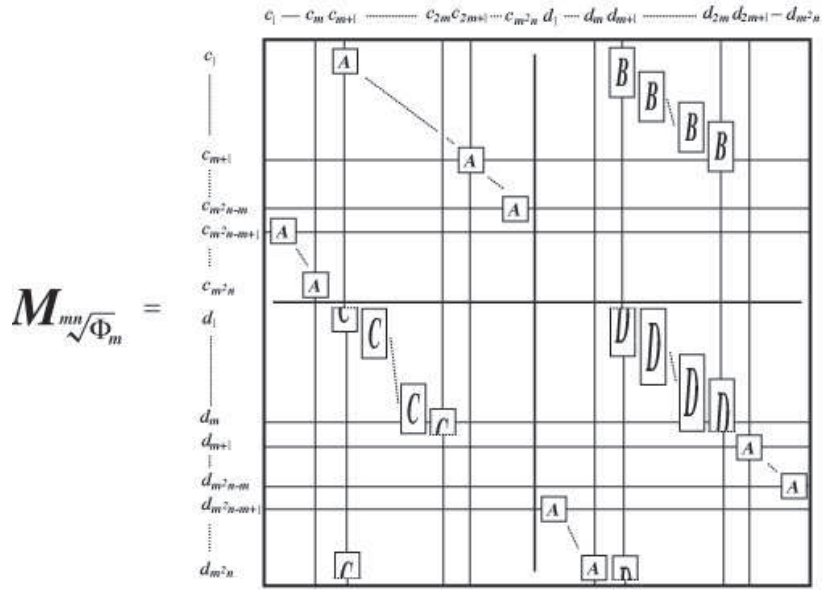, width=5.5in} \caption{The curve matrix
for $\sqrt[mn]\Phi_{m}$. }
\end{figure}  

Denote by $\lambda_{m}$ 
the Perron root of $M_{\Psi_{m}}$ and by 
\[ \upsilon_{m} \;\;=\;\; (a_{1}\dots ,a_{m^{2}n},b_{1},\dots , b_{m^{2}n})\] 
the corresponding
Perron vector, normalized to have $L^{1}$ norm 1 {\em i.e.}\ so that $\upsilon_{m}$
is a probability vector.  
Let $\upsilon_{m}^{\rm avg}$ be the vector 
\[  \upsilon_{m}^{\rm avg}\;\; =\;\; (a_{1}+\dots +a_{m^{2}n},b_{1}+\dots +b_{m^{2}n}) ,\]
which is also a probability vector.  In the
case of $\Phi$ and $\tilde{\Phi}$, the Perron roots are identical and will
be denoted $\lambda$; if $\upsilon = (a,b)^{\sf T}$ is the $L^{1}$ norm 1 Perron vector
of $M_{\Phi}$, $\tilde{\upsilon}=(1/2m^{2}n)(a,\dots , a,b,\dots ,b)^{\sf T}$
is the $L^{1}$ norm 1 Perron vector for $M_{\tilde{\Phi}}$.  
We recall that this spectral data
has the following interpretation:

\begin{enumerate}
\item  The Perron roots $\lambda_{m}$, $\lambda$ are equal to
the entropies of $\Psi_{m}$, $\Phi$.
\item  Let $\tau_{m}$ be the bigon track formed from the curve families 
$\mathcal{C}$, $\mathcal{D}$.
The measures $\mu_{m}$, $\tilde{\mu}$ formed from the Perron vectors
$\upsilon_{m}$, $\tilde{\upsilon}$, parametrize the unstable
laminations $\mathfrak{f}^{u}_{m}$,
$\tilde{\mathfrak{f}}^{u}$ in the cone ${\sf C}(\tau_{m} )$.
\end{enumerate}

Note that the column sums of $M_{(\Psi_{m})^{n}}$ have uniform
upper and lower bounds $B$ and $b>1$.  We thus obtain the bound \cite{Ga}
\[ 1<b<(\lambda_{m})^{n} <B .\]  
We may then assume, after passing to a subsequence if necessary, that
the $\lambda_{m}$ converge to some value $\lambda^{\ast}>1$.
We shall need to control the following entries
of $\upsilon_{m}$:

\begin{clai}\label{entriesgoto0}  $a_{m},a_{m+1},b_{m},b_{m+1}\rightarrow 0$ as $m\rightarrow\infty$.
\end{clai}

\noindent {\it Proof of claim.}  Let $\xi_{m}$ be the Perron root of 
$M_{\sqrt[nm]{\Phi}_{m}}$:
thus $(\xi_{m})^{m}=\lambda_{m}$.
If one of the four entries listed in the statement
does not converge to 0, then consideration of the
matrix $M_{(\Psi_{m})^{n}}$ shows that none of them do.
Thus, let us suppose that $a_{m}\not\rightarrow 0$.  It follows
that eventually $a_{2m}\geq \delta $ for some positive $\delta$.
Examination of the action of $M_{\sqrt[nm]{\Phi}_{m}}$ on $\upsilon_{m}$
shows that $a_{m+im}=\xi_{m}^{i-1}a_{2m}$
for $i=2,\dots mn$.  However $\xi_{m}^{i-1}>1$ for all $i$,
and since $\upsilon_{m}$ is a probability vector, this would imply that 
$mn\delta<mna_{2m}<1$, impossible since $m\rightarrow\infty$.  This proves
the claim.

\vspace{2mm}

Let us shorten notation, 
writing $\hat{\mathfrak{f}}=\hat{\mathfrak{f}}^{u}$ and
$\hat{\mathfrak{f}}_{m}=\hat{\mathfrak{f}}_{m}^{u}$ for
the unstable foliations of $\hat{\Phi}$ and $\hat{\Psi}_{m}$.  
Let $\hat{c}$ be any simple closed solenoid in $\mathcal{C}_{Z}(\widehat{\Sigma})$.

Then
\[
\text{\bf I}\big( (\lambda_{m})^{n}\hat{\mathfrak{f}}_{m},\hat{c}\big)\;\;=\;\;
\text{\bf I}\big( (\Psi_{m})^{n}\hat{\mathfrak{f}}_{m},\hat{c}\big)\;\;=\;\;
\int_{\widehat{E}} \widehat{M_{(\Psi_{m})^{n}}\upsilon_{m}} |\hat{\tau}\cap\hat{c}| ,
 \]
 where $\widehat{M_{(\Psi_{m})^{n}}\upsilon_{m}}$ is the lift of the vector
 $M_{(\Psi_{m})^{n}}\upsilon_{m}$ to $\widehat{E}$.  Now since $\hat{c}$ is the
 lift of a simple closed curve $c$ in $Z$, we have that
\[
\text{\bf I}\big( \hat{\mathfrak{f}}_{m},\hat{c}\big)
\;\;=\;\; \frac{1}{\deg (Z)}\cdot
\text{\bf I}\big( \mathfrak{f}^{\rm avg}_{m},c\big) ,
\] 
where $\mathfrak{f}^{\rm avg}_{m}$ is the measured lamination in $Z$
corresponding to the weight $\upsilon_{m}^{\rm avg}$.  However an
examination of the matrices $M_{\tilde{\Phi}}$ and $M_{(\Psi_{m})^{n}}$
yields
\[ (\lambda_{m})^{n}\upsilon_{m}^{\rm avg}\;\; =\;\; M_{\Phi} \upsilon_{m}^{\rm avg} +
\epsilon_{m}\] where
\[\epsilon_{m}\;\; =\;\; \big((rN)^{2}(a_{m}+a_{m+1}),-(rN)^{2}b_{m+1}+
(rN)^{3}a_{m}+
((rN)^{2}-(rN)+1)a_{m+1}+(rN)^{4}b_{m}\big)^{T}.\] 
By Claim~\ref{entriesgoto0}, it follows that $\epsilon_{m}\rightarrow (0,0)^{T}$
or $\upsilon_{m}^{\rm avg}$ converges to an eigenvector of $M_{\Phi}$ of eigenvalue
$\lambda_{\ast}^{n}$.  But since $\lambda_{\ast}^{n}> 1$ and the second
eigenvalue of $M_{\Phi}$ is strictly less than 1, we must have that $\upsilon_{m}^{\rm avg}
\rightarrow\upsilon$ and $\lambda_{\ast}=\lambda^{1/n}$.  In particular,
\[ \lim_{m\rightarrow\infty}
\text{\bf I}\big( \mathfrak{f}^{\rm avg}_{m},c\big)
\;\; =\;\; \text{\bf I}\big( \mathfrak{f},c\big) .\]
This takes care of the unstable part of the theorem; the stable part
is handled by repeating the above argument for $\Phi^{-1}$ and $\Psi_{m}^{-1}$.
\end{proof}

\begin{theo}  Every pseudo Anosov $\hat{\Phi}$ has an $L^{1}$ $n$th root
for all $n$.
\end{theo}

\begin{proof}  Since $\mathcal{P}$ is directionally dense,
there exists a sequence 
$\{ \hat{\Phi}^{(g)} \} \subset\mathcal{P}$,
where $\hat{\Phi}^{(g)}$ is the lift of a pseudo Anosov
$\Phi^{(g)}:X_{g}\rightarrow X_{g}$ in which $X_{g}$ is a surface of genus 
$g\rightarrow\infty$, and 
the axes $A_{g}\rightarrow A$ = axis of $\hat{\Phi}$.  For each $g$, let 
$\{ \hat{\Psi}^{(g)}_{m}\}$
be the sequence of pseudo Anosovs constructed above.  Then by Lemma~\ref{Zintconv} 
we may obtain an $L^{1}$ $n$th root $\{ \hat{\Psi}_{m}\}$
of $\hat{\Phi}$ by extracting a suitable diagonal subsequence of 
$\{ \hat{\Psi}^{(g)}_{m}\}$.  Indeed, a suitable diagonal subsequence 
$\{ \hat{\Psi}_{m}\}$ yields
a sequence of pseudo Anosov homeomorphisms whose stable and unstable laminations
intersection converge to those of $\hat{\Phi}$.  By Theorems~\ref{HubbMas} and
\ref{intconvimpl1}, this gives rise to a sequence of quadratic differentials 
$\hat{q}_{i}$
along the associated axes $A_{i}$ which $L^{1}$-converge to the quadratic
differential $\hat{q}$ determined by $\hat{\mathfrak{f}}^{u}$ and 
$\hat{\mathfrak{f}}^{s}$.  
\end{proof}

\bibliographystyle{amsalpha}

\end{document}